\documentstyle[12pt]{article}
\newtheorem{theorem}{Theorem}[section]
\newtheorem{corollary}[theorem]{Corollary}
\newtheorem{definition}{Definition}
\newtheorem{example}[theorem]{Example}

\newtheorem{lemma}[theorem]{Lemma}
\newtheorem{proposition}[theorem]{Proposition}
\title{Some covering properties of the $\alpha$-topology\thanks
{AMS (1991) Subject Classification --- Primary: 54A05, 54D20;
Secondary: 54D30, 54H05. \newline Key words and phrases:
$\alpha$-set, $S$-closed, para-$S$-closed, $rc$-Lindel\"{o}f,
sg-compact, $\alpha$-(sub)paracompact. \newline Research
supported partially by the Ella and Georg Ehrnrooth Foundation
at Merita Bank, Finland.}}
\author{Francisco G. Arenas, Jiling Cao, Julian 
Dontchev and Maria Luz Puertas}
\date{}
\begin{document}
\baselineskip=20pt
\newcommand{\fxy}{$f \colon (X,\tau) \rightarrow (Y,\sigma)$}
\newcommand{\xt}{\mbox{$(X,{\cal T})$}}
\newcommand{\xta}{\mbox{$(X,{\cal T}^{\alpha})$}}
\maketitle
\begin{abstract}
Recently, Mr\v{s}evi\'{c} and Reilly discussed some covering
properties of a topological space and its associated
$\alpha$-topology in both topological and bitopological ways. The
main aim of this paper is to investigate some common and
controversial covering properties of $\cal T$ and ${\cal
T}^{\alpha}$.
\end{abstract}

\section{Introduction}\label{s1}

In 1965, Nj{\aa}stad introduced the notion of an $\alpha$-set in
a topological space $(X, \cal T)$, and proved that the collection
of all $\alpha$-sets in $(X, \cal T)$ is a topology on $X$, finer
that $\cal T$. Having two related topologies on the same
underlying set, it is quite natural to ask whether they share
some topological properties. The sharing of separation axioms has
been considered by Dontchev \cite{Do}, Jankovi\'{c} and Reilly
\cite{JR}, Mr\v{s}evi\'{c} and Reilly \cite{MR1}, etc.\ as a part
of their investigations. Recently, Mr\v{s}evi\'{c} and Reilly
\cite{MR2} discussed some covering properties of $(X, \cal T)$
and $(X, {\cal T}^{\alpha})$ in both topological and
bitopological ways. As is shown in \cite{MR2}, $\cal T$ and
${\cal T}^{\alpha}$ do not share compactness, countable
compactness, Lindel\"{o}fness and paracompactness. On the other
hand, Cao and Reilly \cite{CR} proved that $(X, \cal T)$ is
almost compact (resp.\ almost paracompact) if and only if $(X,
{\cal T}^{\alpha})$ is almost compact (resp.\ almost
paracompact). In light of these results, we study the
$S$-closedness of $(X, {\cal T}^{\alpha})$ and show in
Section~\ref{s2} that $\cal T$ and ${\cal T}^{\alpha}$ share most
of $S$-closed-like properties. Some controversial properties are
discussed in Section~\ref{s3}.

Let $(X, \cal T)$ be a topological space and $A \subseteq X$. The
interior (resp.\ closure) of $A$ is denoted by ${\rm Int} A$
(resp.\ ${\rm Cl} A$). Recall that $A$ is {\em semi-open} (resp.\
{\em regular open, an $\alpha$-set}) if $A \subseteq {\rm
Cl}({\rm Int} A)$ (resp.\ $A = {\rm Int}({\rm Cl} A)$, $A
\subseteq {\rm Int}({\rm Cl}({\rm Int} A))$). In addition, $A$
is called {\em semi-closed\/} (resp.\ {\em regular closed,
$\alpha$-closed\/}) if its complement $X \setminus A$ is
semi-open (resp.\ regular open, an $\alpha$-set). The {\em
semi-closure\/} (resp.\ {\em $\alpha$-closure\/}) of $A$, denoted
by ${\rm sCl} A$ (resp.\ ${\rm Cl}_{\alpha} A$), is defined as
the intersection of all semi-closed (resp.\ $\alpha$-closed)
subsets containing $A$. The $\alpha$-semi-closure of $A$, denoted
by ${\rm sCl}_{\alpha} A$ is defined in a similar way. The family
of all semi-open sets (resp.\ $\alpha$-sets, regular closed sets)
of $(X, \cal T)$ is denoted by $SO(X, \cal T)$ (resp.\ ${\cal
T}^{\alpha}$, $RC(X, \cal T)$). No separation axioms are assumed
unless it is explicitly stated.

A topological space $(X, {\cal T})$ is {\em semi-compact}
\cite{Dor} if every cover of $X$ by semi-open sets has a finite
subcover. Moreover, $(X, {\cal T})$ is {\em $S$-closed} \cite{Th}
(resp.\ {\em $s$-closed} \cite{DN}) if for every cover of $X$ by
semi-open sets there is a finite subfamily whose closures (resp.\
semi-closures) form a cover of $X$. It is easy to show that $(X,
{\cal T})$ is $S$-closed if and only if every cover by regular
closed sets has a finite subcover. In a similar manner, a
topological space $(X, \cal T)$ is called {\em $rc$-Lindel\"{o}f}
\cite{JK} if every cover of $X$ by regular closed subsets has a
countable subcover.

\begin{lemma}\label{le2.1}
{\em \cite{Nj}} $SO(X, {\cal T}^{\alpha}) = SO(X, {\cal T})$ for
any space $(X, {\cal T})$. $\Box$
\end{lemma}

Semi-open sets are not the only classes of sets that \xt\ and
\xta\ have always in common. The following classes of sets are
also shared by both topologies in question: locally dense (=
preopen) sets \cite{CM1}, nowhere dense sets, dense and codense
sets, clopen sets, $\beta$-open (= semi-preopen) sets \cite{PN1}
and of course $\alpha$-open sets.

Recall that a subset $A$ of a topological space $(X,\cal T)$ is
called {\em sg-open} \cite{BL1} (resp.\ {\em g-open} \cite{L1})
if every semi-closed (resp.\ closed) subset of $A$ is included
in the semi-interior (resp.\ interior) of $A$. Complements of
sg-open (resp.\ g-open) sets are called {\em sg-closed} (resp.\
{\em g-closed}). The family of all sg-closed (resp.\ g-closed)
subsets of a topological space \xt\ is denoted by $SGC\xt$
(resp.\ $GC\xt$). Although the definitions of sg-closed and
g-closed sets are very similar to each other, sg-closed and
g-closed sets behave in a very different way. More precisely, 
sg-closed sets are more close (and related) to semi-closed and
$\beta$-closed sets, while g-closed sets behave more like the
other types of generalized closed sets, i.e.\ gs-closed,
gp-closed, $\delta$-g-closed (see \cite{DM1} for more details).
We will see next that this type of behavior is valid also in
connection with the $\alpha$-topology.

\section{$S$-closed-like properties}\label{s2}

\begin{proposition}\label{p1}
Let $(X, {\cal T})$ be a topological space. Then:

{\rm (1)} ${\rm sCl}_{\alpha} A = {\rm sCl} A$ for any $A
\subseteq X$.

{\rm (2)} $SGC\xt = SGC\xta$.
\end{proposition}

{\em Proof.} (i) follows easily from Lemma~\ref{le2.1}.

(ii) Let $A \in SGC\xt$ and let $A \subseteq U$, where $U$ is
semi-open in \xta. By Lemma~\ref{le2.1}, $U$ is semi-open also
in \xt. Thus, ${\rm sCl} A \subseteq U$. By (1), ${\rm
sCl}_{\alpha} A \subseteq U$. Thus $A \in SGC\xta$. Conversely,
assume that $A \in SGC\xta$ and $A \subseteq U$, where $U$ is
semi-open in \xt. Again by Lemma~\ref{le2.1}, $U$ is semi-open
in \xta. Thus, ${\rm sCl}_{\alpha} A \subseteq U$. Using again
(1), we have ${\rm sCl} A \subseteq U$. This shows that $A \in
SGC\xt$. $\Box$

We provide an example showing that in general $GC\xt \not=
GC\xta$.

\begin{example}\label{e1}
{\em Let $X = \{ a,b,c \}$ and ${\cal T} = \{ \emptyset, \{ a \},
X \}$. Set $A = \{ a,b \}$. Observe that $A$ is g-closed in \xt.
Since ${\cal T}^{\alpha} = \{ \emptyset, \{ a \}, \{ a,b \}, \{
a,c \}, X \}$, then it is easily checked that $A$ is not
g-closed in \xta.}
\end{example}

\begin{lemma}\label{le2.2}
{\em \cite{Ja}} Let $(X, {\cal T})$ be a topological space. Then
${\rm Cl}_{\alpha} A = {\rm Cl} A$ for any $A \in SO(X, {\cal
T})$. $\Box$
\end{lemma}

\begin{proposition}\label{pr2.1}
$RC(X, {\cal T}^{\alpha}) = RC(X, {\cal T})$ for any space $(X,
{\cal T})$.
\end{proposition}

{\em Proof.} Let $A \in RC(X, {\cal T})$. Then $A= {\rm Cl} G$
for some $G \in {\cal T}$. By Lemma~\ref{le2.1}, $A ={\rm
Cl}_{\alpha} G$. Hence $A \in RC(X, {\cal T}^{\alpha})$.
Conversely, let $A \in RC(X, {\cal T}^{\alpha})$. Then  $A= {\rm
Cl}_{\alpha} G$ for some $G \in {\cal T}^{\alpha}$. Since each
$\alpha$-set is semi-open, we have $A= {\rm Cl} G = {\rm Cl}({\rm
Int} G)$. Furthermore, $G \subseteq A \subseteq {\rm Cl} G$
implies ${\rm Cl}({\rm Int} A) = {\rm Cl}({\rm Int} G)$.
Therefore, $A ={\rm Cl}({\rm Int} A)$, and $A \in RC(X, {\cal
T})$. $\Box$

\bigskip

Recently, the concept of a sg-compact space was introduced
independently by Caldas in \cite{CC1} and by Devi, Balachandran
and Maki in \cite{DBM1}. A topological space $(X,{\cal T})$ is
called {\em sg-compact} \cite{CC1} if every cover of $X$ by
sg-open sets has a finite subcover. In \cite{DBM1}, sg-compact
spaces are called {\em $SGO$-compact}. Sg-compact spaces are
studied in detail in \cite{DG2}.

As consequences of Lemma~\ref{le2.1}, Proposition~\ref{p1} (ii)
and Proposition~\ref{pr2.1}, we have the following.

\begin{theorem}\label{th2.1}
Let $(X, {\cal T})$ be a topological space. Then:

(a) $(X, {\cal T}^{\alpha})$ is semi-compact if and only if $(X,
{\cal T})$ is semi-compact.

(b) $(X, {\cal T}^{\alpha})$ is $S$-closed if and only if $(X,
{\cal T})$ is $S$-closed.

(c) $(X, {\cal T}^{\alpha})$ is $s$-closed if and only if $(X,
{\cal T})$ is $s$-closed.

(d) $(X, {\cal T}^{\alpha})$ is $rc$-Lindel\"{o}f if and only if
$(X, {\cal T})$ is $rc$-Lindel\"{o}f.

(e) $(X, {\cal T}^{\alpha})$ is sg-compact if and only if $(X,
{\cal T})$ is sg-compact. $\Box$
\end{theorem}

Recall that a subset $A$ of $(X, {\cal T})$ is {\em $S$-closed
relative to $X$} (resp.\ {\em $s$-closed relative to $X$}) if for
every cover of $A$ by semi-open sets in $(X, \cal T)$, there
exists a finite subfamily whose closures (resp.\ semi-closures)
in $(X, \cal T)$ from a cover of $A$. Furthermore, $(X, \cal T)$
is {\em locally $S$-closed} \cite{No} (resp.\ {\em locally
$s$-closed} \cite{Ba}) if each point of $X$ has a neighbourhood
which is $S$-closed relative to $X$ (resp.\ $s$-closed relative
to $X$).

\begin{corollary}
Let $(X, {\cal T})$ be a topological space. Then:

(a) $(X, {\cal T}^{\alpha})$ is locally $S$-closed if and only
if $(X, {\cal T})$ is locally $S$-closed.

(b) $(X, {\cal T}^{\alpha})$ is locally $s$-closed if and only
if $(X, {\cal T})$ is locally $s$-closed. $\Box$
\end{corollary}

In \cite{Ch}, Chen defined a space $(X, \cal T)$ to be {\em
para-$S$-closed} if every cover of $X$ by semi-open sets has a
locally finite refinement by semi-open sets whose union is dense
in $X$. In a recent paper \cite{JK}, Jankovi\'{c} and
Konstadilaki introduced the notion of para-$rc$-Lindel\"{o}fness.
Recall that a topological space $(X, \cal T)$ is {\em
para-$rc$-Lindel\"{o}f} if every cover of $X$ by regular closed
sets has a locally countable refinement by regular closed sets.
Next, we prove that $\cal T$ and ${\cal T}^{\alpha}$ share these
properties.

\begin{theorem}\label{th2.2}
Let $(X, {\cal T})$ be a topological space. Then the following
conditions are equivalent:

(a) $(X, {\cal T})$ is para-$S$-closed.

(b) Every cover of regular closed sets of $(X, {\cal T})$ has a
locally finite refinement consisting of regular closed sets of
$(X, {\cal T})$.

(c) Every cover of regular closed sets of $(X, {\cal
T}^{\alpha})$ has a locally finite refinement consisting of
regular closed sets of $(X, {\cal T}^{\alpha})$.

(d) $(X, {\cal T}^{\alpha})$ is para-$S$-closed.
\end{theorem}

{\em Proof.} $(a) \Rightarrow (b)$. Let ${\cal F} = \{
F_{\gamma}: \gamma \in \Delta \}$ be a cover by regular sets of
$(X, {\cal T})$. Then it is also a cover by semi-open sets of
$(X, {\cal T})$. Therefore, it has a locally finite refinement
${\cal V}$ by semi-open sets of $(X, {\cal T})$ such that $X =
{\rm Cl}(\cup {\cal V})$. Then, $\{ {\rm Cl} V: V \in {\cal
V} \}$ is a locally finite refinement of ${\cal F}$ consisting
of regular closed sets of $(X, {\cal T})$.

$(b) \Rightarrow (c)$. It is obvious.

$(c) \Rightarrow (d)$. Let ${\cal U} = \{ U_{\gamma}: \gamma \in
\Delta \}$ be a cover by semi-open sets of $(X, {\cal
T}^{\alpha})$. Then $\{ {\rm Cl} U_{\gamma}: \gamma \in \Delta
\}$ is cover by regular closed sets of $(X, {\cal T}^{\alpha})$.
Thus, it has a locally finite refinement ${\cal F}$ by regular
closed sets of $(X, {\cal T}^{\alpha})$. Without loss of
generality, we may assume that ${\cal F} = \{ F_{\gamma}: \gamma
\in \Delta \}$ such that $F_{\gamma} \subseteq {\rm Cl}
U_{\gamma}$ for each $\gamma \in \Delta$. Set $V_{\gamma} =
F_{\gamma} \cap {\rm Int} U_{\gamma}$ for each $\gamma \in
\Delta$, and ${\cal V}=\{ V_{\gamma}: \gamma \in \Delta \}$. Then
${\cal V} \subseteq SO(X, {\cal T}^{\alpha})$. Moreover,
$F_{\gamma} \subseteq {\rm Cl} V_{\gamma}$ for each $\gamma \in
\Delta$. Hence, ${\cal V}$ is a locally finite refinement of
${\cal U}$ consisting of semi-open sets of $(X, {\cal
T}^{\alpha})$ such that $X = {\rm Cl}(\cup {\cal V})$.

$(d) \Rightarrow (a)$. Let ${\cal U} = \{ U_{\gamma}: \gamma \in
\Delta \}$ be a cover by semi-open sets of $(X, {\cal T})$. Then
${\cal U}$ is also a cover by semi-open sets of $(X, {\cal
T}^{\alpha})$. Thus, it has a locally finite refinement ${\cal
V}$ of semi-open sets of $(X, {\cal T}^{\alpha})$ such that $X
= {\rm Cl}(\cup {\cal V})$. Note that ${\cal V} \subseteq SO(X,
{\cal T})$ and ${\cal V}$ is locally finite in $(X, {\cal T})$.
It follows that $(X, {\cal T})$ is para-$S$-closed. $\Box$

\bigskip

Similar to Theorem~\ref{th2.2}, we can obtain the following.

\begin{theorem}\label{th2.3}
Let $(X, \cal T)$ be a topological space. Then $(X, {\cal
T}^{\alpha})$ is para-$rc$-Lindel\"{o}f if and only if $(X, \cal
T)$ is para-$rc$-Lindel\"{o}f. $\Box$
\end{theorem}

It is well known that every regular Lindel\" of space is
paracompact. Analogous to this, we have the following result.
Recall that a topological space $(X,\tau)$ is called {\em 
extremally disconnected} if the closure of every open subset of
$X$ is also open or equivalently if every regular closed set is
regular open.

\begin{theorem}\label{t29} 
Every extremally disconnected, $rc$-Lindel\" of space is
para-$S$-closed.
\end{theorem}

{\em Proof.} Suppose that $(X, \cal T)$ is an extremally
disconnected and $rc$-Lindel\" of space. Let $\cal U$ be a cover
of $X$ by regular closed sets. Then $\cal U$ has a countable
subcover $\{ U_n: n \in \omega \}$. Define $V_n = U_n \setminus
\cup_{k=1}^{n-1} U_k$ for each $n \in \omega$. Then, it is easy
to see $\{ V_n: n \in \omega \} \subseteq RC(X, \cal T)$ and $V_n
\subseteq U_n$ for each $n \in \omega$. For each $x \in X$, let
$n(x) = min \{ n \in \omega: x \in U_n \}$. Clearly, we have $x
\in V_{n(x)}$. It follows that $\{V_n: n \in \omega \}$ is a
cover of $X$ by regular closed sets. Since $(X, \cal T)$ is
extremally disconnected, $U_{n(x)}$ is an open neighbourhood of
$x$ for each $x \in X$. On the other hand, $U_{n(x)} \cap V_n =
\emptyset$ for all $n > n(x)$. Therefore, $\{V_n: n \in \omega
\}$ is a locally finite refinement of $\cal U$. By
Theorem~\ref{th2.2}, $(X, \cal T)$ is para-$S$-closed. $\Box$

\section{$\alpha$-subparacompact spaces}\label{s3}

In this last section we prove a subspace theorem for
$\alpha$-subparacompact spaces.

\begin{definition}\label{d31}
{\em A topological space \xt\ is called {\em
$\alpha$-subparacompact} if every $\alpha$-open cover of $X$ has
a $\sigma$-discrete closed refinement.}
\end{definition}

Clearly, every $\alpha$-subparacompact space is subparacompact
but not vice versa as the following example shows:

\begin{example}\label{e31}
{\em Let $X$ be the real line with topology in which the only
nontrivial open set is $\{ 0 \}$. Note that $\{ \{ 0,y \} \colon
y \not= 0 \}$ is an $\alpha$-open cover of $X$ which has no
$\sigma$-discrete closed refinement. Thus, even a compact space
need not be $\alpha$-subparacompact.}
\end{example}

Next, we provide an example of a connected, Tychonoff,
$\alpha$-subparacompact space which is not even metacompact.

\begin{example}\label{e32}
{\em Recall that a measurable set $E \subseteq {\bf R}$ has
density $d$ at $x \in {\bf R}$ if $$\lim_{h \rightarrow 0}
\frac{m(E \cap [x-h,x+h])}{2h}$$ exists and is equal to $d$. Set
$\phi(E) = \{ x \in {\bf R} \colon d(x,E) = 1 \}$. The open sets
of the density topology $\tau_d$ are those measurable sets $E$
that satisfy $E \subseteq \phi(E)$. Note that every nowhere dense
subset of the density topology is closed \cite{T1}. Hence every
$\alpha$-open set is open. Thus the subparacompactness of density
topology \cite{T1} implies automatically its
$\alpha$-subparacompactness. On the other hand, the density
topology is not paracompact, in fact it is not even metacompact
\cite{T1}.}
\end{example}

Recall that a subset $A$ of a topological space \xt\ is called
a {\em generalized $\alpha$-closed set} (briefly {\em
g$\alpha$-closed}) \cite{MDB2} if ${\rm Cl}_{\alpha}(A) \subseteq
U$, whenever $A \subseteq U$ and $U$ is $\alpha$-open. We call
a subset $A$ of a topological space $(X,\tau)$
{\em $F_{\sigma}$-$g{\alpha}$-closed} if $A$ is countable union
of g$\alpha$-closed subsets of $X$. The set of all rationals $\bf
Q$ (in the Real line) is an example of an
$F_{\sigma}$-$g{\alpha}$-closed set which is not
g$\alpha$-closed. 

\begin{theorem}\label{t32}
Let $X$ be an $\alpha$-subparacompact space, and let $A$ be
$F_{\sigma}$-$g{\alpha}$-closed. Then $A$ is
$\alpha$-subparacompact (as a subspace), in particular,
$\alpha$-subparacompactness is a $\alpha$-closed hereditarily.
\end{theorem}

{\em Proof.} Let $A = \cup_{n \in \omega} A_n$, where each $A_n$
is g$\alpha$-closed. Let ${\cal U} = \{ U_i \colon i \in I \}$
be a cover of $\alpha$-open subsets of $(A,{\cal T}|A)$. Note
that for each $i \in I$, there exists $V_i \in {\cal T}^{\alpha}$
(i.e.\ $V_i$ is $\alpha$-open in \xt) such that $V_i \cap A =
U_i$. Since union of $\alpha$-open sets is $\alpha$-open, then
$V = \cup_{i \in I} V_i$ is $\alpha$-open in \xt. Since each
$A_n$ is g$\alpha$-closed, then ${\rm Cl}_{\alpha} (A_n)
\subseteq V$. Observe that $\{ X \setminus {\rm Cl}_{\alpha}
(A_n): n \in \omega \} \cup \{ V_i \colon i \in I \}$ is an
$\alpha$-open cover of \xt. Since $X$ is $\alpha$-subparacompact,
then there exists a $\sigma$-discrete closed refinement, say
${\cal W} = \cup_{m \in \omega} {\cal W}_{m}$. For $m, n \in
\omega$, set ${\cal W}_{mn}^{'} = \{ W \cap A_n \colon W \in
{\cal W}_{m} \}$. Clearly, $\cup_{m \in \omega} \cup_{n \in
\omega} {\cal W}_{mn}^{'}$ is a $\sigma$-discrete closed
refinement of $\cal U$ in $(A,{\cal T}|A)$. This shows that $A$
is $\alpha$-subparacompact subspace of \xt. $\Box$

\begin{corollary}
Every closed subspace of an $\alpha$-subparacompact space is also
$\alpha$-subparacompact. $\Box$
\end{corollary}

There is a result due to Burke \cite{Bu1} cited in Theorem 5.2
of \cite{Y1} (see also the reference to that in Theorems
2.14 -- 2.16 (pp.\ 224) of Junnila's survey \cite{Ju1}): A space
is subparacompact if and only if each open cover has a
$\sigma$-closure preserving closed refinement. In a similar
fashion one can prove the following.

\begin{lemma}\label{lfm1}
A topological space $(X,\tau)$ is $\alpha$-subparacompact if and
only if each $\alpha$-open cover has a $\sigma$-closure
preserving closed refinement. $\Box$
\end{lemma}

Recall that a function \fxy\ is called {\em $\alpha$-irresolute}
\cite{MT1} if the preimage of every $\alpha$-open subset of
$(Y,\tau)$ is $\alpha$-open in $(X,\tau)$.

\begin{theorem}\label{fm1}
Every closed $\alpha$-irresolute image of an
$\alpha$-subparacompact space is also $\alpha$-subparacompact.
\end{theorem}

{\em Proof.} Let \fxy\ be a closed $\alpha$-irresolute (not
necessarily continuous) map from the $\alpha$-subparacompact
space $X$ onto the topological space $Y$ and let ${\cal V} = \{
V_i :i \in I \}$ be an $\alpha$-open cover of $Y$. From the
$\alpha$-irresoluteness of $f$, we have that ${\cal U} = \{
f^{-1}(V_i) \colon i \in I \}$ is an $\alpha$-open cover of $X$.
Since $X$ is $\alpha$-subparacompact, then $\cal U$ has a
$\sigma$-closure preserving closed refinement, that is, there
exists ${\cal Z} = \cup_{i=1}^{\infty}{\cal Z}_{i}$, where each
${\cal Z}_i$ is a closure preserving family of closed sets and
$\cal Z$ refines $\cal U$. Since $f$ is closed, then by
Lemma~\ref{lfm1} ${\cal W} = \cup_{i=1}^{\infty}{\cal W}_{i}$,
where ${\cal W}_{i} = \{ f(Z) \colon Z \in {\cal Z}_{i}$, is a
$\sigma$-closure preserving closed family in $Y$. It is
straightforward to check that $\cal W$ is a refinement of $\cal
V$. $\Box$

The product of an $\alpha$-subparacompact space and a compact
space need not to be $\alpha$-subparacompact (take an infinite
indiscrete space that is clearly $\alpha$-subparacompact and the
compact non-$\alpha$-subparacompact space from
Example~\ref{e31}). 

If \xta\ is compact, then \xt\ is usually called {\em
$\alpha$-compact}. Properties of $\alpha$-compact spaces were
studied in 1986 by Noiri and Di Maio \cite{NDM1}.

{\bf Question 1.} Is the product of two $\alpha$-subparacompact
spaces $\alpha$-subparacompact? Is the product of an
$\alpha$-subparacompact space and an $\alpha$-compact space
necessarily $\alpha$-subparacompact?

\begin{definition}
{\em A topological space is called {\em $\alpha$-paracompact} if
every $\alpha$-open cover of $X$ has a locally finite open
refinement.}
\end{definition}

\begin{proposition}
If \xt\ is a Hausdorff $\alpha$-paracompact space, then \xta\ is
normal, in particular if \xt\ is Hausdorff and
$\alpha$-paracompact, then ${\cal T}^{\alpha} = {\cal T}$, i.e.,
$X$ is a nodec space.
\end{proposition}

{\em Proof.} Let $A \subseteq X$ be an $\alpha$-closed set and
$x \not\in A$. For every $y \in A$ there exists an open set $U_y$
such that $y \in U \setminus y$ and $x \not\in \overline{U_y}$.
Then $\{ U_y \colon y \in A \} \cup \{ X \setminus A \}$ is an
$\alpha$-open cover of $X$. Let ${\cal W}$ be a locally finite
open refinement and let $U = \cup \{ W \in {\cal W} \colon W \cap
A \neq \emptyset \}$. Then $U$ is open, contains $A$ and
$\overline{U} = \cup \{ \overline{W} \colon W \cap A \neq
\emptyset \}$. But each such set $W$ is contained in some $U_y$,
and hence $\overline{W} \subseteq \overline{U_y}$ and thus $x
\not\in \overline{U}$. Now let $A,B$ disjoint $\alpha$-open
subsets of $X$. For each $x \in A$ there exists an open set $V_x$
such that $x \in V_x$ and $\overline{V_x} \cap B = \emptyset$.
Then $\{ V_x \colon x \in A \} \cup \{ X \setminus A \}$ is an
$\alpha$-open cover of $X$. Let ${\cal W}$ be a locally finite
open refinement and let $V = \cup \{ W \in {\cal W} \colon W \cap
A \neq \emptyset \}$. Then $V$ is open, contains $A$ and
$\overline{U} = \cup \{ \overline{W} \colon W \cap A \neq
\emptyset \}$. But each such set $W$ is contained in some $V_x$,
and hence $\overline{W} \subseteq \overline{V_x}$ and thus
$\overline{U} \cap B = \emptyset$. $\Box$

\begin{theorem}
Suppose $(X,{\cal T})$ is Hausdorff and $\alpha$-paracompact.
Then $(X,{\cal T}^{\alpha})$ is Hausdorff and paracompact.
\end{theorem}

{\em Proof.} By a theorem of Engelking \cite[Theorem 5.1.5, page
373]{E1}, $(X,{\cal T}^{\alpha})$ is normal. By a result of
Dontchev \cite{JD1}, ${\cal T} = {\cal T}^{\alpha}$. $\Box$

{\bf Question 2.} Let $(X,{\cal T}^{\alpha})$ be subparacompact.
Is $(X,{\cal T})$ $\alpha$-subparacompact?

\
\begin{center}
Area of Geometry and Topology\\Faculty of
Science\\Universidad de Almer\'{\i}a\\04071
Almer\'{\i}a\\Spain\\e-mail: {\tt farenas@ualm.es}
\end{center}
\
\begin{center}
Department of Mathematics\\The University of Auckland\\Private
Bag 92019\\Auckland\\New Zealand\\e-mail: {\tt
cao@math.auckland.ac.nz}
\end{center}
\
\begin{center}
Department of Mathematics\\University of Helsinki\\PL 4,
Yliopistonkatu 15\\00014 Helsinki\\Finland\\e-mail: {\tt
dontchev@cc.helsinki.fi}, {\tt
dontchev@e-math.ams.org}\\http://www.helsinki.fi/\~{}dontchev/
\end{center}
\
\begin{center}
Area of Geometry and Topology\\Faculty of
Science\\Universidad de Almer\'{\i}a\\04071
Almer\'{\i}a\\Spain\\e-mail: {\tt mpuertas@ualm.es}
\end{center}
\
\end{document}